\pdfoutput=1
\RequirePackage{ifpdf}
\ifpdf 
\documentclass[pdftex]{sigma}
\else
\documentclass{sigma}
\fi

\numberwithin{equation}{section}

\DeclareMathOperator{\av}{av}
\DeclareMathOperator{\At}{At}

\begin{document}

\newcommand{\arXivNumber}{1401.2787}

\allowdisplaybreaks

\renewcommand{\PaperNumber}{070}

\FirstPageHeading

\ShortArticleName{On the Conjectures Regarding the 4-Point Atiyah Determinant}

\ArticleName{On the Conjectures Regarding\\ the 4-Point Atiyah Determinant}

\Author{Mazen N.~BOU KHUZAM~$^\dag$ and Michael J.~JOHNSON~$^\ddag$}

\AuthorNameForHeading{M.N.~Bou Khuzam and M.J.~Johnson}

\Address{$^\dag$~American University of Iraq, Suleimaniya, Street 10, Quarter 410,\\
\hphantom{$^\dag$}~Ablakh area Building no.~7 Sul, Iraq}
\EmailD{\href{mailto:mazen.boukhuzam@auis.edu.iq}{mazen.boukhuzam@auis.edu.iq}}

\Address{$^\ddag$~Department of Mathematics, Faculty of Science, Kuwait University, Kuwait}
\EmailD{\href{mailto:yohnson1963@hotmail.com}{yohnson1963@hotmail.com}}

\ArticleDates{Received January 15, 2014, in f\/inal form June 23, 2014; Published online July 05, 2014}

\vspace{-1.5mm}

\Abstract{For the case of 4~points in Euclidean space, we present a~computer aided proof of Conjectures~II and~III made
by Atiyah and Sutclif\/fe regarding Atiyah's determinant along with an elegant factorization of the square of the
imaginary part of Atiyah's determinant.}


\Classification{51K05; 51P99}

\vspace{-2.5mm}

\section{Introduction}

The Atiyah determinant is a~complex-valued determinant function $\At(P_{1},\dots,P_{n})$ associated with~$n$ distinct
points $P_{1},\dots,P_{n}$ of ${\mathbb R}^{3}$.
It was constructed by M.F.~Atiyah in~\cite{A1} in his attempt at answering a~natural geometric question posed
in~\cite{BR} and arising from the study of the spin statistics theorem using classical quantum theory.
The original conjecture of Atiyah was that $\At$ does not vanish for all conf\/igurations of distinct points $P_{1},\dots,
P_{n} \in {\mathbb R}^{3}$.
The conjecture was verif\/ied in the linear case (all points lie on a~straight line) and in the case $n=3$ by Atiyah
in~\cite{A1}.
However, the case $n\geq 4$ turned out to be notoriously dif\/f\/icult.
In a~subsequent paper~\cite{AS}, Atiyah and Sutclif\/fe studied the function $\At$ and added two new conjectures (after
normalizing~$\At$) which imply the original conjecture of Atiyah.
They provided compelling numerical evidence of the validity of all three conjectures.
The three conjectures can be stated as follows: For all distinct points $P_{1},\dots,P_{n}$ of ${\mathbb R}^{3}$ (and all
$n\geq 1$) we have:
\begin{alignat*}{3}
& \text{(I)}\quad &&
\At(P_{1},\dots,P_{n}) \neq 0, &
\\
& \text{(II)}\quad &&
|\At(P_{1},\dots,P_{n})| \geq \prod\limits_{i<j} (2 r_{ij}),
\qquad
\text{where}
\quad
r_{ij}=||\overrightarrow{P_{i}P_{j}}||, &
\\[-1mm]
& \text{(III)}\quad &&
 |\At(P_{1},\dots,P_{n})|^{n-2} \geq \prod\limits_{k=1}^{n} |\At(P_{1},\dots,P_{k-1}, P_{k+1},\dots,P_{n})|.
\end{alignat*}

From the statement of these conjectures we can see that $\text{III} \implies \text{II} \implies \text{I}$.
The three conjectures have been very resistant since their inauguration time.
The f\/irst conjecture was proved by Eastwood and Norbury~\cite{EN} for the case $n=4$.
Other attempts were successful only on special conf\/igurations (see~\cite{D} and~\cite{MP}).
In this paper, we build on the work of Eastwood and Norbury by presenting a~computer aided proof of Conjectures~II and~III in the case $n=4$ and we also give an elegant factorization of the square of the imaginary part of the Atiyah determinant.

The construction of the determinant is as follows: One starts with~$n$ distinct points $P_{1},\dots, P_{n}$ $\in
{\mathbb R}^{3}$.
By considering $P_{j}$ as an observer of the other $n-1$ points we obtain $n-1$ vectors
$\overrightarrow{P_{j}P_{1}},\dots,\overrightarrow{P_{j}P}_{j-1}$,
$\overrightarrow{P_{j}P}_{j+1},\dots,\overrightarrow{P_{j}P_{n}}$ in ${\mathbb R}^{3}$.
We lift each of these vectors from ${\mathbb R}^{3}$ to ${\mathbb C}^{2}$ using the Hopf map $h: {\mathbb C}^{2} \rightarrow
{\mathbb R}^{3}$ given by $h(z,w)=((|z|^{2}-|w|^{2})/2, z\overline{w})$ to obtain $n-1$ points of ${\mathbb C}^{2}$.
Note that the lifts are not unique and are def\/ined up to phase because $h(\lambda z, \lambda w)= |\lambda |^{2} h(z,w)$.
Consequently, our lifts can be considered as points of ${\mathbb C}P^{1}$.
Taking the symmetric product of these lifts gives a~vector $V_{j}$ in ${\mathbb C}P^{n}$ because $\odot_{n}
{\mathbb C}P^{1}={\mathbb C}P^{n}$.
Atiyah's f\/irst conjecture was that $\{V_{1},\dots, V_{n} \}$ is a~linearly independent set.
In other words, the determinant of the matrix having the vector $V_{j}$ as its $jth$ column is nonzero.
This determinant is well-def\/ined up to a~phase factor.
To get rid of the phase factor ambiguity, we apply the following normalization imposed by Atiyah: If $(z,w)$ is the
chosen lift of $\overrightarrow{P_{i}P_{j}}$ and $i<j$, then $(-\overline{w}, \overline{z})$ must be the lift of
$\overrightarrow{P_{j}P_{i}}$.
After this normalization, this determinant is called the Atiyah determinant and is denoted by~$\At$.

It is immediate from the above construction that~$\At$ is coordinate free and is independent of solid motion.
In other words, the determinant function~$\At$ is invariant under translations and rotations in ${\mathbb R}^{3}$.
Furthermore, the Atiyah determinant is built so that it is independent of the order of the points.
In other words, if $(j_{1},\dots, j_{n})$ is a~permutation of $(1,\dots, n)$ then $\At(P_{j_{1}},\dots,
P_{j_{n}})=\At(P_{1},\dots, P_{n})$.
Another property is that $\At$ gets conjugated under a~plane ref\/lection of the points (see~\cite{A1}).
As a~consequence, $\At$ must be real-valued if the set of points $\{P_{1},\dots, P_{n}\}$ is symmetric relative to a~plane
(e.g.\
if the points are co-planar) since a~ref\/lection in the plane leaves the set of points unchanged.

Let us start computing $\At$ in the cases $n=2$ and $n=3$.
For the case $n=2$, we have two distinct points~$A$ and~$B$.
We can identify ${\mathbb R}^{3}$ with ${\mathbb R}\times {\mathbb C}$ and assume (possibly after a~solid motion) that~$A$ and~$B$
have coordinates $(0,0)$ and $(0,x)$ respectively, where $x>0$ is the distance from~$A$ to~$B$.
By choosing $(\sqrt{x}, \sqrt{x})$ as a~lift of $\overrightarrow{AB}$, we are forced to take $(-\sqrt{x}, \sqrt{x})$ as
a~lift of $\overrightarrow{BA}$.
Consequently, Atiyah's determinant is:
\begin{gather*}
\At(A,B)= \left|
\begin{matrix}
\sqrt{x} & -\sqrt{x}
\\
\sqrt{x} & \sqrt{x}
\end{matrix}
\right| =2x,
\qquad
\text{where}
\quad
x=||\overrightarrow{AB}||.
\end{gather*}

\begin{figure}[t]
\centering
\includegraphics{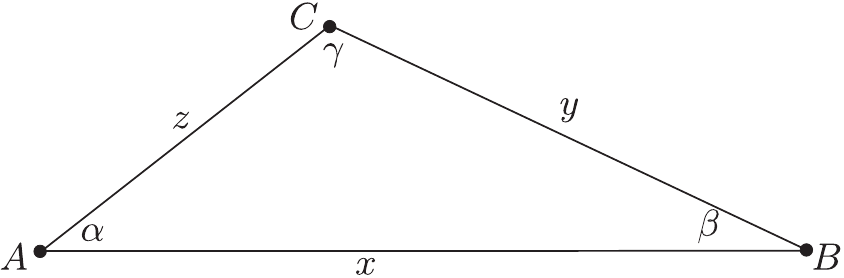}
\caption{Three points.}\label{Fig1}
\end{figure}

Let us now consider the case $n=3$.
Assume (possibly after a~solid motion) that $A=(0,0)$, $B=(0,x)$, and $C=(0,z e^{I\alpha})$ where~$I$ denotes
$\sqrt{-1}$, $y=||\overrightarrow{BC}||$, $z=||\overrightarrow{AC}||$, $x=||\overrightarrow{AB}||$
and~$\alpha$,~$\beta$,~$\gamma$ are the angles indicated in Fig.~\ref{Fig1}.
When the f\/irst point is considered as a~vision point we obtain $\overrightarrow{AB}=(0,x)$, $\overrightarrow{AC}=(0,z e^{I\alpha})$
whose lifts under the Hopf map~$h$ are $(\sqrt{x}, \sqrt{x})$ and $(\sqrt{z}, \sqrt{z} e^{-I \alpha})$.
And when $B=(0,x)$ is the vision point, we obtain the vectors $\overrightarrow{BA}=(0,-x)$ and
$\overrightarrow{BC}=(0,-y e^{-I \beta})$ whose lifts are $(-\sqrt{x}, \sqrt{x})$ and $(-\sqrt{y}, \sqrt{y} e^{I\beta})$.
Similarly, the lifts corresponding to the vision point~$C$ are $(-\sqrt{z} e^{I \alpha},\sqrt{z})$ and $(-\sqrt{y} e^{-I\beta},-\sqrt{y})$.
The symmetric tensor product of the vectors are then $\sqrt{xz} (1, 1+e^{-I\alpha},  e^{-I\alpha})$,
$\sqrt{xy} (1, -1-e^{I\beta},  e^{I\beta})$ and $\sqrt{yz} (e^{I (\alpha -\beta)}, e^{I \alpha}-e^{-I \beta},-1)$, respectively.
Consequently, we obtain the Atiyah determinant for three points as
\begin{gather*}
\At(A,B,C)=xyz \left|
\begin{matrix}
1 & 1 & e^{I (\alpha -\beta)}
\\
1+e^{-I\alpha} & -1-e^{I\beta} & e^{I \alpha}-e^{-I \beta}
\\
e^{-I\alpha} & e^{I\beta} & -1
\end{matrix}
\right|
\end{gather*}

This determinant expands to $xyz[6+2(\cos \alpha +\cos \beta +\cos \gamma)]$, which can be written as $xyz[8+8 \sin
\frac{\alpha}{2} \sin \frac{\beta}{2} \sin \frac{\gamma}{2}]$.
Using the identity $\sin \frac{\alpha}{2}=\frac{1}{2} \sqrt{\frac{(a+b-c)(a+c-b)}{bc}}$ and similar identities for $\sin
\frac{\beta}{2}$ and $\sin \frac{\gamma}{2}$, we can rewrite the Atiyah determinant for three points as
\begin{gather}
\label{eq1.1} \At(A,B,C)=8xyz+d_3(x,y,z),
\end{gather}
where $d_3$ is the polynomial def\/ined by $d_3(x,y,z)=(-x+y+z)(x+y-z)(x+y-z)$.
From the triangle inequality it follows that $d_3(x,y,z)$ is nonnegative, and so Conjecture III is verif\/ied for three
points.
(Note that there is no need to use $|\At|$ because in this case $\At$ is real.)

\section{The case of four points}

\begin{figure}[t]
\centering
\includegraphics{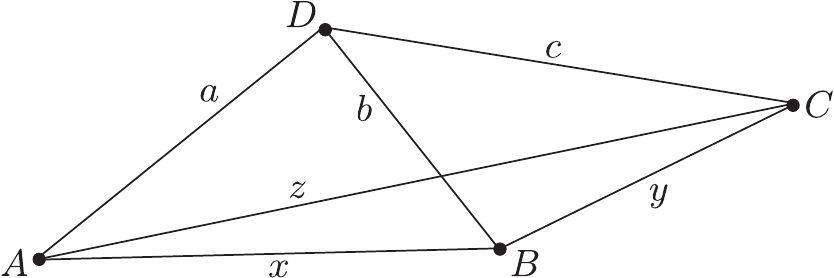}
\caption{Four points.}\label{Fig2}
\end{figure}

Given four points $A$, $B$, $C$, $D$ in ${\mathbb R}^3$, the vector $u=U(A,B,C,D)$ in ${\mathbb R}^6$, called the {\it vector of pair-wise
distances}, is def\/ined by $u=(a,b,c,x,y,z)$ where $a=||\overrightarrow{AD}||$, $b=||\overrightarrow{BD}||$,
$c=||\overrightarrow{CD}||$, $x=||\overrightarrow{AB}||$, $y=||\overrightarrow{BC}||$, $z=||\overrightarrow{AC}||$ (see
Fig.~\ref{Fig2}).
The function~$U$, as def\/ined above, maps ${\mathbb R}^{3\times 4}$ into ${\mathbb R}^6$, and it is clear that~$U$ is neither
injective nor surjective.
A~vector $u\in{\mathbb R}^6$ is said to be {\it geometric} if it belongs to the range of~$U$.
For convenience, we adopt the convention that $\At(A,B,C,D)$ equals $0$ when the points $A$, $B$, $C$, $D$ are not distinct, and
Conjecture~II (for four points) becomes
\begin{gather*}
\text{(II)}
\quad
|\At(A,B,C,D)|\geq 64abcxyz
\quad
\text{for all points}
\quad
A,B,C,D \in {\mathbb R}^3.
\end{gather*}

Atiyah's determinant is designed to be invariant under permutations of the points.
Each of the 24 possible permutations of the four points $A$, $B$, $C$, $D$ results in a~permutation of the pair-wise distances
$a$, $b$, $c$, $x$, $y$, $z$.
Specif\/ically, if $u=(a,b,c,x,y,z)\in{\mathbb R}^6$, then the 24 resultant permutations are
\begin{alignat*}{5}
& u_0=(a,b,c,x,y,z), \ && u_1=(a,x,z,b,y,c), \ && u_2=(b,c,a,y,z,x), \ && u_3=(x,b,y,a,c,z),&
\\
& u_4=(c,a,b,z,x,y), \ && u_5=(z,y,c,x,b,a), \ && u_6=(y,z,c,x,a,b), \ && u_7=(c,b,a,y,x,z),&
\\
& u_8=(x,y,b,z,c,a), \ && u_9=(a,c,b,z,y,x), \ && u_{10}=(z,x,a,y,b,c), \ && u_{11}=(b,a,c,x,z,y),&
\\
& u_{12}=(z,c,y,a,b,x), \ && u_{13}=(x,z,a,y,c,b), \ && u_{14}=(x,a,z,b,c,y), \ && u_{15}=(y,x,b,z,a,c),&
\\
&u_{16}=(y,b,x,c,a,z), \ && u_{17}=(y,c,z,b,a,x), \ && u_{18}=(c,y,z,b,x,a), \ && u_{19}=(z,a,x,c,b,y),&
\\
& u_{20}=(b,x,y,a,z,c), \ \  && u_{21}=(c,z,y,a,x,b), \ \ && u_{22}=(a,z,x,c,y,b), \ \ && u_{23}=(b,y,x,c,z,a).&
\end{alignat*}

A function $f:{\mathbb R}^6 \to {\mathbb R}$ is said to be {\it symmetric} if $f(u)=f(u_i)$ for $i=0,1,\ldots, 23$ and is {\it
skew-symmetric} if $f(u)=(-1)^i f(u_i)$ for $i=0,1,\ldots, 23$.
The {\it symmetric average} of~$f$ is the symmetric function $\av[f]$ def\/ined by
\begin{gather*}
\av[f](u) = \frac1{24}\sum\limits_{i=0}^{23} f(u_i).
\end{gather*}

Using Maple, Eastwood and Norbury have found that the real part of $\At(A,B,C,D)$ can be expressed
as $\Re \At(A,B,C,D) = d_4(u)$, where $d_4$ is the homogeneous polynomial of degree $6$ given~by
\begin{gather}
\label{eq2.1} d_4(u)=60p_4(u) + 4 n_4(u) + 2 z_4(u) + 12\av\big[a\big((b+c)^2-y^2\big)d_3(x,y,z)\big],
\end{gather}
where $p_4(u)=abcxyz$, $d_3$ is def\/ined in \eqref{eq1.1}, $n_4(u)=p_4(u)- d_3(xc,ay,bz)$ and
\begin{gather*}
z_4(u)=a^2y^2\big(b^2 + c^2 +x^2 + z^2\big)+b^2z^2\big(a^2+c^2+x^2+y^2\big)+c^2x^2\big(a^2+b^2+y^2+z^2\big)
\\
\phantom{z_4(u)=}
{}-\big(a^4y^2 + a^2 y^4 + b^4 z^2 + b^2 z^4 + c^4 x^2 + c^2 x^4\big)
\\
\phantom{z_4(u)=}
{}-\big(a^2 b^2 x^2 + a^2 c^2 z^2 + b^2 c^2 y^2 + x^2 y^2 z^2\big).
\end{gather*}
Eastwood and Norbury use the notation $144V^2$ in place of $z_4(u)$.
If $u=U(A,B,C,D)$, the value $z_4(u)$ equals $144V^2$, where~$V$ denotes the volume of the tetrahedron formed by the
points $A$, $B$, $C$, $D$, and it therefore follows that $z_4(u)\geq 0$.
It would be erroneous to infer from this that the polynomial $z_4$ is nonnegative on all of ${\mathbb R}^6$; the above
statement implies only that $z_4$ is nonnegative on geometric vectors.

Having expressed $\Re \At(A,B,C,D)=d_4(u)$ as in~\eqref{eq2.1}, Eastwood and Norbury then invoke the inequalities $z_4(u)\geq 0$,
$(b+c)^2\geq y^2$, $d_3(x,y,z)\geq0$ and $abcxyz\geq d_3(xc,ay,bz)$ (i.e.\ $n_4(u)\geq 0$) to conclude that
\begin{gather*}
|\At(A,B,C,D)|\geq \Re \At(A,B,C,D)=d_4(u) \geq 60p_4(u),
\end{gather*}
which proves Conjecture~I and comes close to proving Conjecture~II.

Regarding the imaginary part of $\At(A,B,C,D)$, Eastwood and Norbury have shown that its square can be written as $(\Im
\At(A,B,C,D))^2=F_4(u)$, where $F_4$ is a~symmetric homogeneous polynomial of degree $12$.
Whereas $d_4$ seems unwilling to be expressed in a~simple manner, we have found that $F_4$ factors elegantly as
\begin{gather*}
F_4 = w_4^2 z_4,
\end{gather*}
where $w_4$ is the skew-symmetric homogeneous polynomial of degree $3$ given by
\begin{gather*}
w_4(u) = \big(a^2+y^2\big)(b-c-x+z) + \big(b^2+z^2\big)(-a+c+x-y) + \big(c^2+x^2\big)(a-b+y-z)
\\
\phantom{w_4(u) =}
{}+2(cx+yz)(-a+b)+2(ay+xz)(-b+c)+2(bz+xy)(a-c).
\hspace{1.3 cm} 
\end{gather*}
Note that since $w_4$ is skew-symmetric it follows that $w_4^2$ is symmetric.
As mentioned in the introduction, the imaginary part of $\At(A,B,C,D)$ vanishes whenever the set of four points
$\{A,B,C,D\}$ is symmetric about a~plane.
Interestingly, this property can be derived from the above factorization: Assuming $\{A,B,C,D\}$ is symmetric about
a~plane, it then follows that $u=U(A,B,C,D)$ is invariant under some odd permutation of the four points $A$, $B$, $C$, $D$.
Since $w_4$ is skew-symmetric, we must have $w_4(u)=0$ and hence $F_4(u)=0$.
\section{A linear program related to Conjecture II}
Since $|\At(A,B,C,D)|\geq \Re \At(A,B,C,D)=d_4(u)$, in order to prove Conjecture II, it suf\/f\/ices to show that the
polynomial $d_4$ satisf\/ies
\begin{gather}
\label{eq3.1}
d_4(u)\geq 64p_4(u)
\quad
\text{for all geometric vectors}
\quad
u.
\end{gather}

If one has in hand a~collection $f_1,f_2,\ldots,f_k$ of symmetric homogeneous polynomials of degree $6$ which are known
to be nonnegative on geometric vectors, then one can `have a~go' at \eqref{eq3.1} by solving the linear program
\begin{gather}
\begin{split}
& \text{maximize}
\quad
\alpha,
\\
& \text{subject to}
\quad
d_4 = \alpha p_4 + \sum\limits_{j=1}^k \lambda_j f_j,
\quad
\text{with}
\quad
\lambda_1,\lambda_2,\ldots,\lambda_k\geq 0.
\end{split}
\label{eq3.2}
\end{gather}
If \eqref{eq3.2} is feasible and if the optimal objective value is $\alpha=64$ (we will see later that $\alpha>64$ is
impossible), then we immediately obtain \eqref{eq3.1}.
The remaining dif\/f\/iculty is that of f\/inding suitable polynomials $\{f_j\}$.
One means of generating a~large collection of such polynomials, which we now describe, stems from the triangle
inequality.

The four points $A$, $B$, $C$, $D$ contain four (possibly degenerate) triangles and each triangle, by means of the triangle
inequality, gives rise to three linear polynomials which are nonnegative when $u=(a,b,c,x,y,z)$ is geometric.
For example, the triangle $A$, $B$, $C$ yields $-x+y+z$, $x-y+z$ and $x+y-z$.
In all, there are twelve such linear polynomials which we refer to as {\it triangular variables} and use the notation
$t=(t_1,t_2,\ldots,t_{12})$, where
\begin{alignat*}{5}
& t_1 = -a+b+x, \qquad && t_4 = -b + c +y, \qquad && t_7 = -a +c +z,\qquad && t_{10} = -x + y +z,&
\\
& t_2 = a-b+x, \qquad && t_5 = b -c +y, \qquad && t_8 = a -c +z, \qquad && t_{11} = x -y +z,&
\\
& t_3 = a+b-x, \qquad && t_6 = b + c -y, \qquad && t_9 = a + c -z, \qquad && t_{12} = x + y -z.&
\end{alignat*}
A~vector \looseness=-1 $\alpha\in{\mathbb Z}^{12}_+$ is called a~{\it multi-index} with {\it order} $|\alpha|=\alpha_1 +\alpha_2 + \dots +
\alpha_{12}$.
Employing the standard notation $t^\alpha=t_1^{\alpha_1} t_2^{\alpha_2} \cdots t_{12}^{\alpha_{12}} $, we see that
$t^\alpha$ represents a~homogeneous polynomial of degree $|\alpha|$ in the variables $(a,b,c,x,y,z)$.
Applying the symmetric average, we conclude that $\av[t^\alpha]$ represents a~symmetric homogeneous polynomial of
degree $|\alpha|$ which is nonnegative on geometric vectors.
For integers $\ell\geq 0$, we def\/ine ${\mathbb T}_\ell$ to be the set of all polynomials $\av[t^\alpha]$ with
$|\alpha|=\ell$:
\begin{gather*}
{\mathbb T}_\ell=\{\av[t^\alpha]:|\alpha|=\ell\}.
\end{gather*}

Numerically, we have found that if one chooses $\{f_j\}$ equal to ${\mathbb T}_6$, then the linear program~\eqref{eq3.2} is feasible
and has optimal objective value $\alpha=32$.
The formulation \eqref{eq2.1} of Eastwood and Norbury can be understood in the context of~\eqref{eq3.2} as the result of including, in
addition to ${\mathbb T}_6$, the two symmetric polynomials $z_4$ and $n_4$ which are nonnegative on geometric vectors.
Numerically solving~\eqref{eq3.2} with $\{f_j\}$ equal to $\{z_4,n_4\}\cup {\mathbb T}_6$, we have found that the optimal objective
value is $\alpha=60$, and~\eqref{eq2.1} is indeed an optimal solution of~\eqref{eq3.2} as the term $\av[a((b+c)^2-y^2)d_3(x,y,z)]$
can be written as a~nonnegative linear combination of polynomials in~${\mathbb T}_6$.

In order to further increase the optimal objective value~$\alpha$ in \eqref{eq3.2}, we need other symmetric polynomials which
are nonnegative on geometric vectors.
In pursuit of this, we have identif\/ied the following twenty-one geometric vectors~$u$ where $d_4(u)=64p_4(u)$ (all are
obtained as $u=U(A,B,C,D)$ with $A$, $B$, $C$, $D$ collinear or non-distinct):
\begin{alignat}{4}
& (0, 1, 4, 1, 4, 4), \qquad && (0, 4, 8, 4, 7, 8), \qquad && (0, 6, 0, 6, 6, 0),&
\nonumber\\
& (0, 1, 1, 1, 2, 1), \qquad && (0, 5, 5, 5, 5, 5), \qquad && (0, 8, 8, 8, 1, 8),&
\nonumber\\
& (0, 1, 3, 1, 4, 3), \qquad && (0, 6, 3, 6, 8, 3), \qquad && (0, 6, 7, 6, 3, 7),&
\nonumber\\
& (0, 6, 6, 6, 9, 6), \qquad && (0, 1, 1, 1, 0, 1), \qquad && (0, 5, 3, 5, 3, 3),& \label{eq3.4}\\
& (3, 3, 1, 0, 2, 2), \qquad && (9, 9, 7, 0, 2, 2),\qquad && (13, 13, 7, 0, 6, 6),&
\nonumber\\
& (19, 11, 7, 8, 4, 12), \qquad && (17, 13, 4, 4, 9, 13), \qquad && (15, 8, 7, 7, 1, 8),&
\nonumber\\
& (9, 8, 1, 1, 7, 8), \qquad && (11, 9, 8, 2, 1, 3), \qquad && (17, 9, 2, 8, 7, 15).&\nonumber
\end{alignat}
Both $d_4$ and $p_4$ vanish on the f\/irst f\/ifteen of these vectors (counting horizontally), but are nonzero on the
remaining six.
In particular, since $d_4(9, 8, 1, 1, 7, 8)=64p_4(9, 8, 1, 1, 7, 8) = 258048 >0$, it follows that there are no feasible
solutions of \eqref{eq3.2} with $\alpha>64$.
On the other hand, if a~feasible solution of \eqref{eq3.2} has been obtained with $\alpha=64$, then it follows that $f_j$
vanishes on all of the vectors in \eqref{eq3.4}, whenever $\lambda_j>0$.
It has been verif\/ied that $z_4$ vanishes on all of these vectors, but $n_4$ does not.
Therefore, the coef\/f\/icient of $n_4$ will be $0$ if~\eqref{eq3.2} has been solved with $\alpha=64$.
We have considered numerous symmetric homogeneous polynomials of degree $6$ which vanish on the vectors in~\eqref{eq3.4}, but
only one of these has resulted in an improvement.
Let $v_4$ denote the skew-symmetric homogeneous polynomial of degree $3$ def\/ined~by
\begin{gather*}
v_4(u)=(b+z-c-x)(c+x-a-y)(a+y-b-z).
\end{gather*}
Then $v_4$ vanishes on the vectors in \eqref{eq3.4}, and numerically solving \eqref{eq3.2} with $\{f_j\}$ equal to
$\{z_4,n_4,v_4^2\}\cup {\mathbb T}_6$, we have found that the optimal objective value is $\alpha=188/3$.
Our obtained identity, which has been verif\/ied in Maple\footnote{The sources of our codes are available at \url{http://www.emis.de/journals/SIGMA/2014/070/codes.zip}.}, is the following:
\begin{gather*}
d_4(u) = \frac{188}3 p_4(u) + \frac{10}3 z_4(u) + \frac43 n_4(u) + \frac23 v_4^2(u) +
\frac13\sum\limits_{|\alpha|=6}\lambda_\alpha \av[t^\alpha],
\end{gather*}
where the six nonzero coef\/f\/icients $\lambda_\alpha$ and corresponding multi-indices~$\alpha$ are given~by
\begin{alignat*}{7}
& \alpha \qquad && \lambda_\alpha \qquad && \alpha \qquad && \lambda_\alpha \qquad && \alpha \qquad && \lambda_\alpha &
\\
& 0 0 0, 0 0 1, 0 1 0, 1 1 2 \qquad && 6 \qquad && 0 0 0, 0 0 1, 0 1 1, 1 1 1 \qquad && 18 \qquad && 0 0 0, 0 0 1, 1 1 0, 1 0 2 \qquad && 6 &
\\
& 0 0 1, 0 0 1, 0 0 1, 1 1 1 \qquad && 14 \qquad && 0 0 1, 0 0 1, 0 1 0, 1 1 1 \qquad && 24 \qquad && 0 0 1, 0 1 1, 1 0 0, 1 1 0 \qquad && 24 &
\end{alignat*}

\section{Proof of Conjecture II for four points}

Let $m_4$ be the symmetric homogeneous polynomial of degree $6$ def\/ined by $m_4 = d_4 - (64p_4 + 4z_4 + v_4^2)$, so that
\begin{gather}
\label{eq4.1}
d_4 = 64p_4 + 4z_4 + v_4^2 + m_4.
\end{gather}
We will show that $m_4$ is nonnegative on geometric vectors, but unfortunately, we have been unable to formulate a~proof
using only polynomials of degree~$6$.
Rather, we have had to multiply $m_4$ by $p_4$ and then work with polynomials of degree $12$.
\begin{theorem}
The product $p_4 m_4$ can be written as a~nonnegative linear combination of polynomials in ${\mathbb T}^{12}$.
\end{theorem}

\begin{proof}
Using Maple, we have verif\/ied that $64 p_4 m_4$ can be written as
\begin{gather}
\label{eq4.2}
64 p_4(u) m_4(u) = \sum\limits_{|\alpha|=12}\lambda_\alpha \av[t^{\alpha}],
\end{gather}
where the sixty-four nonzero coef\/f\/icients $\{\lambda_\alpha\}$ are all positive integers as given in the following
table:
\begin{alignat*}{7}
& \alpha \qquad && \lambda_\alpha \qquad && 0 1 1, 0 2 1, 2 0 1, 1 1 2 \qquad && 6 \qquad && 0 1 1, 1 1 1, 2 2 0, 1 0 2 \qquad && 6 &
\\
& 0 0 1, 0 1 2, 2 1 1, 1 1 2 \qquad && 12 \qquad && 0 1 1, 0 2 1, 2 0 1, 1 2 1 \qquad && 6 \qquad && 0 1 1, 1 1 2, 1 0 2, 2 1 0 \qquad && 12 &
\\
& 0 0 1, 0 1 2, 2 1 1, 1 2 1 \qquad && 12 \qquad && 0 1 1, 0 2 1, 2 0 1, 2 1 1 \qquad && 18 \qquad && 0 1 1, 1 1 2, 1 2 0, 2 1 0 \qquad && 39 &
\\
& 0 0 1, 1 1 2, 1 1 2, 0 1 2 \qquad && 12 \qquad && 0 1 1, 0 2 1, 2 2 1, 1 1 0 \qquad && 54 \qquad && 0 1 1, 1 1 2, 2 1 0, 2 1 0 \qquad && 42 &
\\
& 0 0 1, 1 2 1, 1 2 1, 0 1 2 \qquad && 21 \qquad && 0 1 1, 0 2 2, 0 1 1, 1 2 1 \qquad && 6 \qquad && 0 1 1, 1 2 0, 0 1 1, 2 1 2 \qquad && 27 &
\\
& 0 0 2, 1 1 0, 1 1 1, 2 1 2 \qquad && 9 \qquad && 0 1 1, 0 2 2, 1 1 1, 1 0 2 \qquad && 84 \qquad && 0 1 1, 1 2 0, 0 1 1, 2 2 1 \qquad && 24&
\\
& 0 0 2, 1 1 0, 1 1 1, 2 2 1 \qquad && 9 \qquad && 0 1 1, 0 2 2, 1 1 1, 2 1 0 \qquad && 18 \qquad && 0 1 1, 1 2 0, 0 1 2, 1 1 2 \qquad && 3&
\\
& 0 1 1, 0 1 1, 0 2 1, 2 2 1 \qquad && 30 \qquad && 0 1 1, 0 2 2, 2 1 1, 1 1 0 \qquad && 6 \qquad && 0 1 1, 1 2 1, 0 2 1, 2 0 1 \qquad && 24&
\\
& 0 1 1, 0 1 1, 1 0 1, 2 2 2 \qquad && 56 \qquad && 0 1 1, 1 0 1, 0 2 1, 2 1 2 \qquad && 54 \qquad && 0 1 1, 1 2 1, 1 0 2, 2 1 0 \qquad && 6&
\\
& 0 1 1, 0 1 1, 1 1 0, 2 2 2 \qquad && 48 \qquad && 0 1 1, 1 0 2, 0 1 2, 2 1 1 \qquad && 6 \qquad && 0 1 1, 1 2 1, 1 2 0, 2 1 0 \qquad && 24 &
\\
& 0 1 1, 0 1 1, 1 1 2, 2 2 0 \qquad && 6 \qquad && 0 1 1, 1 0 2, 0 2 2, 1 1 1 \qquad && 36 \qquad && 0 1 1, 2 0 1, 0 1 2, 2 1 1 \qquad && 3&
\\
& 0 1 1, 0 1 1, 1 2 0, 2 1 2 \qquad && 18 \qquad && 0 1 1, 1 0 2, 1 1 2, 1 0 2 \qquad && 18 \qquad && 0 1 1, 2 0 1, 0 2 1, 2 1 1 \qquad && 45&
\\
& 0 1 1, 0 1 1, 1 2 2, 1 0 2 \qquad && 57 \qquad && 0 1 1, 1 0 2, 1 1 2, 2 0 1 \qquad && 24 \qquad && 0 1 2, 0 1 2, 0 1 2, 1 1 1 \qquad && 24 &
\\
& 0 1 1, 0 1 1, 2 0 1, 2 2 1 \qquad && 36 \qquad && 0 1 1, 1 0 2, 2 0 1, 1 2 1 \qquad && 33 \qquad && 0 1 2, 0 1 2, 1 0 2, 1 1 1 \qquad && 8 &
\\
& 0 1 1, 0 1 1, 2 2 2, 0 1 1 \qquad && 6 \qquad && 0 1 1, 1 0 2, 2 1 0, 1 1 2 \qquad && 75 \qquad && 0 1 2, 1 1 1, 0 1 2, 0 2 1 \qquad && 36 &
\\
& 0 1 1, 0 1 2, 1 0 2, 1 1 2 \qquad && 120 \qquad && 0 1 1, 1 0 2, 2 1 0, 1 2 1 \qquad && 12 \qquad && 0 1 1, 0 2 1, 2 2 0, 2 1 0 \qquad && 6&
\\
& 0 1 1, 0 1 2, 1 1 2, 2 1 0 \qquad && 18 \qquad && 0 1 1, 1 1 0, 0 2 1, 2 2 1 \qquad && 48 \qquad && 0 1 1, 1 2 0, 0 1 2, 2 0 2 \qquad && 21&
\\
& 0 1 1, 0 1 2, 1 2 0, 1 2 1 \qquad && 12 \qquad && 0 1 1, 1 1 1, 0 1 2, 2 0 2 \qquad && 21 \qquad && 0 1 1, 1 2 0, 0 2 1, 2 0 2 \qquad && 24&
\\
& 0 1 1, 0 1 2, 2 1 1, 1 0 2 \qquad && 84 \qquad && 0 1 1, 1 1 1, 0 2 1, 2 0 2 \qquad && 18 \qquad && 0 1 1, 2 0 1, 0 1 2, 2 0 2 \qquad && 18&
\\
& 0 1 1, 0 1 2, 2 1 1, 2 0 1 \qquad && 72 \qquad && 0 1 1, 1 1 1, 0 2 1, 2 2 0 \qquad && 54 \qquad && 0 1 1, 2 1 0, 0 2 1, 2 0 2 \qquad && 3&
\\
& 0 1 1, 0 2 1, 0 1 1, 1 2 2 \qquad && 3 \qquad && 0 1 1, 1 1 1, 0 2 2, 1 0 2 \qquad && 75 \qquad && 0 1 2, 0 1 2, 1 2 0, 1 0 2 \qquad && 3&
\\
& 0 1 1, 0 2 1, 1 1 2, 0 1 2 \qquad && 69 \qquad && 0 1 1, 1 1 1, 2 0 2, 2 1 0 \qquad && 12 \qquad && && &\tag*{\qed}
\end{alignat*}
\renewcommand{\qed}{}
\end{proof}
\begin{corollary}
The polynomial $m_4$ is nonnegative on geometric vectors and consequently \eqref{eq3.1} holds, which proves Conjecture~{\rm II} for
four points.
\end{corollary}
\begin{proof}
Let $u=U(A,B,C,D)$ be a~geometric vector.
It follows from \eqref{eq4.2} that $p_4(u) m_4(u)\geq 0$.
If the points $A$, $B$, $C$, $D$ are distinct, then $p_4(u)>0$ and hence $m_4(u)\geq 0$.
On the other hand, if $A$, $B$, $C$, $D$ are not distinct, then they can be approximated by distinct points $A'$, $B'$, $C'$, $D'$ and it
will then follow from the continuity of $m_4$ that $m_4(u)\geq 0$.
\end{proof}

\section{Proof of Conjecture III for four points}

Let $P_4$ denote the symmetric homogeneous polynomial of degree $12$ given~by
\begin{gather*}
P_4(u):=(8xyz+d_3(x,y,z))(8abx+d_3(a,b,x))(8acz+d_3(a,c,z))(8bcy+d_3(b,c,y)),
\end{gather*}
whereby $\At(A,B,C)\At(A,B,D)\At(A,C,D)\At(B,C,D)=P_4(u)$ when $u=U(A,B,C,D)$.
Since $|\At(A,B,C,D)|^2 \geq (\Re \At(A,B,C,D))^2 = d_4^2(u)$, in order to prove Conjecture III, it suf\/f\/ices to show that
\begin{gather}
\label{eq5.1}
d_4^2(u) \geq P_4(u)
\quad
\text{for all geometric vectors}
\quad
u.
\end{gather}

Recall from \eqref{eq4.1} that $d_4$ has been written as $d_4 = 64p_4 + m_4 + (4z_4+v_4^2)$, so it follows that
\begin{gather*}
d_4^2=\big(4z_4+v_4^2\big)d_4 + (64p_4+m_4)d_4
\\
\phantom{d_4^2}
 = \big(4z_4+v_4^2\big)d_4 + (64p_4+m_4)^2 + (64p_4+m_4)\big(4z_4+v_4^2\big)
\\
\phantom{d_4^2}
 = \big(4z_4+v_4^2\big)(d_4 + 32p_4 + m_4) + (64p_4+m_4)^2 + 32p_4\big(4z_4+v_4^2\big).
\end{gather*}
With $M_4$ denoting the symmetric homogeneous polynomial of degree $12$ def\/ined by $M_4 = (64p_4+m_4)^2 +
32p_4(4z_4+v_4^2) - P_4$, we then have
\begin{gather}
\label{eq5.2} d_4^2 = P_4 + \big(4z_4+v_4^2\big)(d_4 + 32p_4 + m_4) + M_4.
\end{gather}

\begin{theorem}
The polynomial $M_4$ is nonnegative on geometric vectors, and consequently \eqref{eq5.1} holds, which proves Conjecture~{\rm III} for
four points.
\end{theorem}
\begin{proof}
Using Maple, we have verif\/ied that $128 M_4$ can be written as
\begin{gather}
\label{eq5.3}
128 M_4(u) = \big(4z_4(u)+v_4^2(u)\big)\sum\limits_{|\alpha| = 6}\mu_\alpha \av[t^\alpha]+\sum\limits_{|\alpha|=12}\nu_\alpha \av[t^\alpha],
\end{gather}
where the coef\/f\/icients $\{\mu_\alpha\}$ and $\{\nu_\alpha\}$ are nonnegative integers: The 6 nonzero coef\/f\/icients
$\mu_\alpha$ and corresponding monomials~$\alpha$ are given in the following table:
\begin{alignat*}{7}
& 0 0 0, 0 0 1, 1 1 1, 1 1 0 \qquad && 1236 \qquad && 0 0 0, 1 0 1, 1 0 1, 1 0 1 \qquad && 3594 \qquad && 0 0 1, 0 1 0, 0 1 1, 0 1 1 \qquad && 300&
\\
& 0 0 0, 1 0 0, 1 0 1, 1 1 1 \qquad && 60 \qquad && 0 0 0, 1 0 1, 1 1 0, 1 0 1 \qquad && 114 \qquad && 0 0 1, 0 1 1, 1 0 1, 0 0 1 \qquad && 1014&
\end{alignat*}

The 114 nonzero coef\/f\/icients $\nu_\alpha$ and corresponding monomials~$\alpha$ are given in the following table:
\begin{alignat*}{7}
& 0 0 0, 1 1 2, 1 2 1, 1 1 2 \qquad && 2019 \qquad && 0 1 1, 1 2 0, 0 2 1, 1 1 2 \qquad && 2184 \qquad && 0 0 1, 1 2 0, 2 0 2, 1 2 1 \qquad && 76 &
\\
& 0 0 1, 0 1 2, 1 2 1, 1 1 2 \qquad && 369 \qquad && 0 1 1, 1 2 1, 0 2 1, 2 0 1 \qquad && 228 \qquad && 0 0 1, 1 2 1, 2 0 1, 2 2 0 \qquad && 6 &
\\
& 0 0 1, 0 1 2, 2 1 1, 1 2 1 \qquad && 138 \qquad && 0 1 1, 1 2 1, 2 1 0, 2 1 0 \qquad && 72 \qquad && 0 0 1, 1 2 1, 2 2 0, 0 2 1 \qquad && 3174&
\\
& 0 0 1, 0 1 2, 2 1 1, 2 1 1\qquad &&666 \qquad && 0 1 1, 2 0 1, 0 2 1, 2 1 1\qquad &&936 \qquad && 0 0 1, 1 2 1, 2 2 0, 1 2 0\qquad &&1266&
\\
& 0 0 1, 0 2 1, 2 1 1, 1 2 1\qquad &&3087 \qquad && 0 1 2, 0 1 2, 0 1 2, 1 1 1\qquad &&3072 \qquad && 0 0 1, 1 2 2, 2 1 0, 1 2 0\qquad &&1428&
\\
& 0 0 1, 0 2 2, 1 1 1, 1 1 2\qquad &&3009 \qquad && 0 1 2, 1 1 1, 0 2 1, 2 0 1\qquad &&1308 \qquad && 0 0 1, 2 0 1, 2 2 0, 1 1 2\qquad &&822&
\\
& 0 0 1, 0 2 2, 1 1 1, 1 2 1\qquad &&1074 \qquad && 0 1 2, 1 1 1, 2 1 0, 0 1 2\qquad &&5374 \qquad && 0 0 1, 2 1 0, 2 2 1, 0 2 1\qquad &&612&
\\
& 0 0 1, 0 2 2, 1 1 2, 1 1 1\qquad &&42 \qquad && 0 0 0, 0 1 2, 1 1 1, 2 2 2\qquad &&60 \qquad && 0 0 1, 2 1 1, 2 2 0, 0 1 2\qquad &&300&
\\
& 0 0 1, 0 2 2, 2 1 1, 1 1 1\qquad &&12114 \qquad && 0 0 0, 0 1 2, 1 1 2, 2 1 2\qquad &&1776 \qquad && 0 0 1, 2 1 1, 2 2 0, 1 2 0\qquad &&3072&
\\
& 0 0 1, 1 0 1, 2 1 1, 1 2 2\qquad &&240 \qquad && 0 0 1, 0 0 1, 1 1 2, 2 2 2\qquad &&138 \qquad && 0 0 2, 0 0 2, 1 1 2, 1 1 2\qquad &&1536&
\\
& 0 0 1, 1 1 1, 1 2 2, 1 0 2\qquad &&4056 \qquad && 0 0 1, 0 0 1, 1 2 2, 2 1 2\qquad &&1398 \qquad && 0 0 2, 0 1 1, 2 2 1, 0 2 1\qquad &&1176&
\\
& 0 0 1, 1 1 1, 2 1 0, 2 2 1\qquad &&444 \qquad && 0 0 1, 0 0 2, 1 2 2, 1 1 2\qquad &&3072 \qquad && 0 0 2, 0 1 2, 1 1 0, 2 1 2\qquad &&1662&
\\
& 0 0 1, 1 2 0, 1 2 1, 1 1 2\qquad &&144 \qquad && 0 0 1, 0 0 2, 2 1 1, 1 2 2\qquad &&1236 \qquad && 0 0 2, 0 1 2, 2 1 0, 1 2 1\qquad &&5136&
\\
& 0 0 1, 1 2 1, 2 1 0, 1 1 2\qquad &&714 \qquad && 0 0 1, 0 1 0, 1 1 2, 2 2 2\qquad &&768 \qquad && 0 0 2, 0 1 2, 2 1 0, 2 1 1\qquad &&762&
\\
& 0 0 1, 1 2 1, 2 1 0, 2 1 1\qquad &&1146 \qquad && 0 0 1, 0 1 0, 1 2 2, 2 1 2\qquad &&384 \qquad && 0 0 2, 0 1 2, 2 1 1, 0 2 1\qquad &&2178&
\\
& 0 0 1, 2 1 1, 2 2 0, 1 1 1\qquad &&1866 \qquad && 0 0 1, 0 1 0, 2 1 2, 1 2 2\qquad &&384 \qquad && 0 0 2, 0 2 2, 1 1 0, 1 1 2\qquad &&1188&
\\
& 0 0 2, 0 1 1, 1 1 1, 1 2 2\qquad &&2808 \qquad && 0 0 1, 0 1 1, 1 0 2, 2 2 2\qquad &&2568 \qquad && 0 0 2, 0 2 2, 1 1 0, 1 2 1\qquad &&150&
\\
& 0 0 2, 0 1 1, 2 1 1, 1 2 1\qquad &&3207 \qquad && 0 0 1, 0 1 1, 2 2 2, 0 1 2\qquad &&1224 \qquad && 0 0 2, 0 2 2, 2 1 1, 0 1 1\qquad &&2178&
\\
& 0 0 2, 0 1 2, 1 1 1, 1 1 2\qquad &&3654 \qquad && 0 0 1, 0 1 2, 1 1 0, 2 2 2\qquad &&2634 \qquad && 0 0 2, 1 0 1, 1 0 1, 2 2 2\qquad &&1245&
\\
& 0 0 2, 0 1 2, 1 1 1, 1 2 1\qquad &&3252 \qquad && 0 0 1, 0 1 2, 1 2 0, 2 1 2\qquad &&66 \qquad && 0 0 2, 1 0 1, 1 2 2, 1 0 2\qquad &&246&
\\
& 0 0 2, 0 2 1, 1 1 1, 1 1 2\qquad &&276 \qquad && 0 0 1, 0 1 2, 1 2 2, 1 0 2\qquad &&822 \qquad && 0 0 2, 1 0 1, 2 2 1, 0 1 2\qquad &&528&
\\
& 0 0 2, 0 2 2, 1 1 1, 1 1 1\qquad &&720 \qquad && 0 0 1, 0 2 1, 1 1 2, 2 2 0\qquad &&840 \qquad && 0 0 2, 1 0 2, 1 2 1, 1 0 2\qquad &&222&
\\
& 0 0 2, 1 1 1, 2 1 2, 0 1 1\qquad &&516 \qquad && 0 0 1, 0 2 2, 1 2 0, 1 1 2\qquad &&4542 \qquad && 0 0 2, 1 1 0, 1 2 2, 1 0 2\qquad &&7548&
\\
& 0 0 2, 1 1 2, 2 1 0, 1 1 1\qquad &&1206 \qquad && 0 0 1, 0 2 2, 2 1 1, 1 2 0\qquad &&2928 \qquad && 0 0 2, 1 2 1, 2 1 0, 0 1 2\qquad &&1074&
\\
& 0 0 2, 1 1 2, 2 1 1, 0 1 1\qquad &&1662 \qquad && 0 0 1, 0 2 2, 2 2 0, 1 1 1\qquad &&6546 \qquad && 0 0 2, 1 2 2, 2 1 0, 0 1 1\qquad &&60&
\\
& 0 1 1, 0 1 1, 0 1 2, 2 1 2\qquad &&10110 \qquad && 0 0 1, 1 0 1, 1 0 2, 2 2 2\qquad &&2898 \qquad && 0 1 1, 0 1 2, 1 2 0, 2 2 0\qquad &&3072&
\\
& 0 1 1, 0 1 1, 0 2 1, 2 1 2\qquad &&1164 \qquad && 0 0 1, 1 0 1, 2 0 1, 2 2 2\qquad &&804 \qquad && 0 1 1, 0 2 1, 2 1 0, 2 2 0\qquad &&342&
\\
& 0 1 1, 0 1 1, 2 0 1, 2 1 2\qquad &&5472 \qquad && 0 0 1, 1 0 1, 2 2 0, 1 2 2\qquad &&1398 \qquad && 0 1 1, 0 2 1, 2 2 0, 1 2 0\qquad &&4668&
\\
& 0 1 1, 0 1 2, 0 1 2, 2 1 1\qquad &&342 \qquad && 0 0 1, 1 0 2, 1 1 0, 2 2 2\qquad &&690 \qquad && 0 1 1, 0 2 1, 2 2 0, 2 1 0\qquad &&4608&
\\
& 0 1 1, 0 1 2, 1 1 2, 2 0 1\qquad &&2178 \qquad && 0 0 1, 1 0 2, 1 2 2, 1 0 2\qquad &&1770 \qquad && 0 1 1, 0 2 2, 1 2 0, 1 0 2\qquad &&768&
\\
& 0 1 1, 0 1 2, 1 2 2, 1 0 1\qquad &&192 \qquad && 0 0 1, 1 0 2, 2 0 2, 1 1 2\qquad &&5532 \qquad && 0 1 1, 0 2 2, 1 2 0, 2 1 0\qquad &&1536&
\\
& 0 1 1, 0 2 1, 2 1 1, 1 2 0\qquad &&5472 \qquad && 0 0 1, 1 0 2, 2 0 2, 1 2 1\qquad &&612 \qquad && 0 1 1, 0 2 2, 2 0 1, 1 2 0\qquad &&1890&
\\
& 0 1 1, 0 2 2, 1 1 1, 1 2 0\qquad &&2376 \qquad && 0 0 1, 1 0 2, 2 1 0, 2 1 2\qquad &&390 \qquad && 0 1 1, 0 2 2, 2 1 0, 1 0 2\qquad &&1152&
\\
& 0 1 1, 0 2 2, 1 2 1, 1 1 0\qquad &&696 \qquad && 0 0 1, 1 0 2, 2 2 1, 1 0 2\qquad &&2634 \qquad && 0 1 1, 0 2 2, 2 1 0, 2 1 0\qquad &&6648&
\\
& 0 1 1, 1 0 1, 0 2 2, 1 1 2\qquad &&2628 \qquad && 0 0 1, 1 1 0, 2 0 1, 2 2 2\qquad &&192 \qquad && 0 1 1, 1 2 0, 0 2 2, 1 0 2\qquad &&10752&
\\
& 0 1 1, 1 0 2, 2 0 1, 1 1 2\qquad &&372 \qquad && 0 0 1, 1 1 2, 2 2 0, 1 2 0\qquad &&600 \qquad && 0 1 1, 2 1 0, 0 1 2, 2 0 2\qquad &&168&
\\
& 0 1 1, 1 1 0, 0 1 2, 1 2 2\qquad &&3558 \qquad && 0 0 1, 1 1 2, 2 2 1, 0 0 2\qquad &&60 \qquad && 0 1 1, 2 1 0, 1 2 0, 2 0 2\qquad &&522&
\\
& 0 1 1, 1 1 2, 0 2 1, 2 0 1\qquad &&948 \qquad && 0 0 1, 1 2 0, 2 0 1, 2 2 1\qquad &&774 \qquad && 0 1 2, 1 2 0, 0 1 2, 2 0 1\qquad &&4440&
\end{alignat*}

It now follows from \eqref{eq5.3} that $M_4$ is nonnegative on geometric vectors and we obtain \eqref{eq5.1} as a~consequence of \eqref{eq5.2}.
\end{proof}

\pdfbookmark[1]{References}{ref}
\LastPageEnding

\end{document}